\documentclass[10pt]{article}
\usepackage{amssymb,amsmath,graphicx,verbatim,eufrak,amsthm,hyperref}
\newtheorem{thm}{Theorem}[section]

\newtheorem{prop}[thm]{Proposition}
\theoremstyle{definition}

\theoremstyle{remark}

\numberwithin{equation}{section}

\newcommand{\ignore}[1]{}
\newcommand{\NB}[1]{{\it Noam: #1}}

\newcommand{\AY}[1]{{\it Ariel: #1}}

\newcommand{\norm}[1]{\left\Vert#1\right\Vert}
\newcommand{\abs}[1]{\left\vert#1\right\vert}

\newcommand{\set}[1]{\left\{#1\right\}}
\newcommand{\setb}[2]{ \left\{#1 \ \Big| \ #2 \right\} }

\newcommand{\Bracket}[1]{\left[#1\right]}
\newcommand{\Soger}[1]{\left(#1\right)}

\newcommand{\Set}[2]{ \left\{#1 \ \big| \ #2 \right\} }
\newcommand{\br}[1]{\left[#1\right]}

\newcommand{\sr}[1]{\left(#1\right)}

\newcommand{\Integer}{\mathbb{Z}}

\newcommand{\Z}{\mathbb{Z}}

\newcommand{\eps}{\varepsilon}

\DeclareMathOperator*{\E}{\mathbb{E}}

\renewcommand{\Pr}{}
\let\Pr\relax
\DeclareMathOperator*{\Pr}{\mathbb{P}}

\def\squareforqed{\hbox{\rlap{$\sqcap$}$\sqcup$}}
\def\qed{\ifmmode\squareforqed\else{\unskip\nobreak\hfil
\penalty50\hskip1em\null\nobreak\hfil\squareforqed
\parfillskip=0pt\finalhyphendemerits=0\endgraf}\fi}

\renewcommand{\th}{^{\mathrm{th}}}
\newcommand{\p}{\partial}

\renewcommand{\P}{\mathcal{P}}

\newcommand{\Ee}{\mathcal{E}}

\begin{document}

\title{Long Range Percolation Mixing Time}

\author{Itai Benjamini\thanks{\texttt{itai.benjamini@weizmann.ac.il}}
\and Noam Berger\thanks{\texttt{berger@math.ucla.edu}} %
\and Ariel Yadin\thanks{\texttt{ariel.yadin@weizmann.ac.il}} }
\date{}

\maketitle

% ---------------------------------------------------------------

\begin{abstract}
We provide an estimate, sharp  up to poly-logarithmic factors, of
the asymptotic almost sure mixing time of the graph created by
long-range percolation on the cycle of length $N$ ($\Integer /
N\Integer$). While it is known that the asymptotic almost sure
diameter drops from linear to poly-logarithmic as the exponent $s$
decreases below $2$ \cite{LRP, Biskup}, the asymptotic almost sure
mixing time drops from $N^2$ only to $N^{s-1}$ (up to
poly-logarithmic factors).
\end{abstract}

% ----------------------------------------------------------------

\section{Introduction}
In this note we study mixing time of simple random walks on the
random graph obtained by adding edges to a graph of a cycle, where
edges are added between any two vertices with probability decaying
with their distance, and independently for any two vertices (in
fact, this probability is $1 - \exp ( -\beta \norm{x-y}^{-s} )$,
where $\beta$ and the exponent $s$ are parameters). See below for
definitions of the model and of mixing times, as well as further
remarks.

The bounds on the mixing time are presented in the coming
sections.

Let us mention two interesting findings. First, the mixing time
undergoes a phase transition as the exponent $s$ decreases below
$2$. Second, in this natural model the almost sure diameter is
a.a.s. (asymptotically almost surely) poly-logarithmic, for
certain range of the parameters, yet the a.a.s. mixing time is
polynomial. Such a gap between the diameter and the mixing time
cannot exist for vertex transitive graphs or in graphs where the
isoperimetric dimension is determined by the volume growth
function.

\ignore{ \NB{Here I would say a.a.s., namely asymptotically almost
surely rather that a.s.} }

The long range percolation graphs gained some recent interest, see
\cite{Kleinberg}. From an algorithmic viewpoint it is useful to
consider the  mixing times of these graphs.

\subsection{The Model}

The model we discuss is the finite long-range percolation model
with polynomial decay. Let $N$ be a positive integer, let $\beta>
0, 1 < s<2$ , and consider the following random graph:  Start with
the cycle on $N$ vertices ($\Integer /N \Integer$). Define
$\norm{x-y} = \min \set{ \abs{x - y} , N - \abs{x - y} }$, (which
is the regular graph-theoretical distance). The following edges
are randomly added:  If $x \neq y$, then $x$ and $y$ will be
attached with probability $1 - \exp ( -\beta \norm{x-y}^{-s} )$. 
The different edges are all independent of each other. The
probability of an edge between two (distant enough) vertices is
very close to $\beta \norm{x-y}^{-s}$. We call the graph created
this way $G_{s,\beta}(N)$.

For updated background on long-range percolation see \cite{Biskup,
Biskup2}.

\subsection{Mixing Time} \label{mixing}

%%%%% Changed

Throughout this paper we consider random walks on a random graph.
In order to avoid issues of convergence to the stationary distribution,
we always consider the 
{\em lazy random walk}, 
i.e. the Markov chain on the vertices
of the graph whose transition matrix is
$P(x,y) = \frac{1}{2\deg(x)}$ if $(x,y)$ is an edge in the graph,
and $P(x,x) = \frac{1}{2}$.  
This insures that the random walk is ergodic and converges to the stationary distribution. 

Also, throughout this paper we consider oriented edges.  Thus, if $x$ and $y$ are adjacent 
vertices of some graph, we consider both $(x,y)$ and $(y,x)$ as edges in the edge set.

%%%%% ======

The \emph{variational distance} of the random walk on a graph $G$
is defined
$$ \Delta_x (t) = \frac{1}{2} \sum_{y \in V(G)} \abs{ P^t (x,y) -
\pi(x) } $$ %
where $P^t$ is the $t\th$ power of the transition matrix of the
walk, and $\pi$ is the stationary distribution, i.e. $\pi(x) =
\frac{\deg(x)}{\abs{E(G)}}$ ($\deg(x)$ is the degree of the vertex
$x$, and recall that $E(G)$ is the set of oriented edges). %%%%% Changed 
$\Delta_x (t)$ measures how close the distribution of the
walk starting at $x$, at time $t$, is to the stationary
distribution $\pi$.

\ignore{

\NB{I'd rather define the variation distance in complete generality, then restrict it to the random walk. Suggested paragraph immediately below:}\\
For two probability measures $\mu$ and $\nu$ on the same space $X$ we say that the variational distance between them is $\Delta(\mu,\nu):=\sup\{\mu(A)-\nu(A):A\subseteq X\}$.
When $X$ is discrete,  $$ \Delta(\mu,\nu) = \frac{1}{2} \sum_{y \in V(G)} \abs{ \mu(\{y\}) -
\nu(\{y\}) }. $$
For an aperiodic Markov kernel $P$ on $X$,  $x\in X$ and $t\geq 1$ integer, we define
$\Delta_x (t):=\Delta(P^t(x,\cdot),\pi)$ with $\pi$ being the stationary distribution.
In our case, $X=V(G)$ and $P$ is the transition kernel for the random walk on $G$.

\AY{I don't see how full generality helps the reader.  I prefer
the less general definition.}

}

Define
$$ \tau_x(\eps) = \min \Set{ t }{ \forall \ t' \geq t \ \Delta_x(t')
\leq \eps } . $$ %
Since we are interested in the time it takes the walk to converge
to the stationary distribution, the \emph{mixing time} is defined
$$ \tau(G) = \max_{x \in V(G)} \tau_x(1/4) $$

It is well know that using the second eigenvalue of the transition
matrix, one can bound the mixing time of the graph.  Formally
(though not in full generality), Diaconis and Strook in \cite{DS}
prove that:
\begin{align} \label{eqn:mixing}
\tau(G) \leq \frac{\log(4 |E(G)|)}{ 1 - \lambda_G } ,
\end{align}
%%%%% Changed
where $\lambda_G$ is the second eigenvalue of $G$
(i.e. the eigenvalues of the transition matrix of the random walk on $G$
are $1 > \lambda_G \geq \lambda_3 \geq \cdots \geq \lambda_{|V(G)|}$.
$1-\lambda_G$ is also called the {\em spectral gap} of $G$.)

For more on mixing see \cite{Aldous, Sinclair}.

%%%%% Changed

% Sinclair in \cite{Sinclair} provides a tool for bounding the
% mixing time of a graph, called \emph{canonical paths method}.
% In this method, we consider oriented edges, i.e. every edge
% $\set{x,y}$ in $G$ is considered as two directed edges $(x,y)$ and
% $(y,x)$.  For any pair of vertices $x \neq y$, we need to specify
% a single path $\gamma_{xy}$ from $x$ to $y$. Sinclair shows that
% for any collection of such paths,
% $$ \frac{1}{1 - \lambda_G} \leq \frac{1}{\abs{E(G)}} \max_{e \in
% E(G)} \sum_{\gamma_{xy} \ni e} d(x) d(y) \abs{\gamma_{xy}} $$ %
% where $E(G)$ is now the set of oriented edges, $\abs{\gamma}$ is
% the length of $\gamma$, and the sum is over all paths
% $\gamma_{xy}$ in the collection containing the edge $e$. So to
% bound the mixing time, it is enough to specify a set of paths, one
% for each pair of vertices, while making sure that any path is not
% too long, and that no edge is overloaded with paths.

%%%%% ======

\subsection{Remarks}

\begin{itemize}

    \item In the following sections we will provide upper and lower bounds
    on the mixing time of $G_{s,\beta}(N)$, that match up to
    poly-logarithmic factors. It is interesting to note that using two
    different methods, we obtain matching bounds.

    \item As will be seen in Section \ref{phase}, 
    a phase transition occurs in the mixing time when
    $s$ passes from below $2$ to above $2$.  Two open questions
    regarding the mixing time are:
    \begin{enumerate}
        \item What is the mixing time at $s=2$.  We note that even
        the diameter is not known in this case.

        \item When $s$ drops below $1$, it is shown in \cite{BKPS}
        that the diameter is bounded.  We conjecture that the
        mixing time is constant (independent of $N$) in this case.
    \end{enumerate}

    \item Long-range percolation gives natural examples of graphs with small 
    diameter (poly-logarithmic, see
    \cite{Biskup}) yet large polynomial mixing time.  
Long range percolation is a natural
model of some social networks, in which the probability you know a
person decays with distance. This suggests that while the diameter
of such networks might be small, sampling from such network via
random walk might take long time.

    \item For almost sure expansion of other models of random graphs
    see \cite{Alon, MPS}. Regarding mixing for random walks on other models of random
    graphs see \cite{BKW, FR} and \cite{BM}.

    \item Another question related to mixing on random graphs is
    the following:  Let $\set{G_n}$ be a family of vertex transitive graphs such
    that $\lim_{n \to \infty} \abs{G_n} = \infty$.  Assume that
    the average degree of the giant component $G'_n$ of Bernoulli
    percolation on $G_n$ is uniformly bounded in $n$.  Prove
    $$ \tau(G'_n) \leq O \sr{ \max \set{ \tau(G_n) \ , \ \log^2
    \abs{G_n} } } . $$

    \item
    Consider uniform spanning tree on the long range percolation graph over $\Integer$.
    Our mixing time estimates show that the mixing time of
    long range percolation at $s=3/2$ is like that of a $4$
    dimensional torus.  Since the transition from tree to forest in the uniform spanning
    tree on $\Z^d$ occurs at $d=4$, this suggests that perhaps the
    uniform spanning tree on the long range percolation graph over
    $\Z$, is supported on a tree a.s. iff $s \geq 3/2$. See
    \cite{BLPS} for background.

\end{itemize}

\section{Upper Bound} \label{upper bound}

%%%%% Changed

\subsection{Multicommodity Flow}

Let $P$ be the transition matrix of a reversible Markov chain,
with stationary distribution $\pi$.
Let $V$ be the set of states of the chain, and let $E$ be the set of
oriented edges; i.e.
$$ E = \set{ (x,y) \in V \times V \ : \ P(x,y) > 0 } . $$

For $x,y \in V$ let $\P(x,y)$ be the set of all simple paths from
$x$ to $y$. Let $\P = \cup_{x \neq y \in V} \P(x,y)$.

A \emph{flow} is a function
$f:\P \to [0,1]$ such that for all $x,y \in V$
$$ \sum_{p \in \P(x,y)} f(p) = \pi(x) \pi(y) . $$
The \emph{edge load} of an edge $e \in E$ is defined as
$$ f(e) = \sum_{\substack{ p \in \P \\ p \ni e}} f(p) |p| . $$
The \emph{congestion} of a flow $f$ is defined as
$$ \rho(f) = \max_{(a,b) \in E} \frac{1}{\pi(a) P(a,b)} f((a,b)) . $$

Theorem 5' of \cite{Sinclair} states that if the eigenvalues of
$P$ are $1 > \lambda \geq \lambda_3 \geq \cdots \geq \lambda_n$
(where $n=|V|$), then for any flow $f$, $(1- \lambda)^{-1} \leq \rho(f)$.
Furthermore, Theorem 8 in \cite{Sinclair} shows that if $P$
induces an ergodic Markov chain (i.e. if $\lambda_n > -1$), then
there exists a flow $f^*$ such that
$\rho(f^*) \leq 16 \tau$, where $\tau$ is the mixing time of the chain.
We call $f^*$ the \emph{optimal flow} for $P$.

\subsection{Upper Bound}

We are now ready to prove an upper bound on the mixing time of $G_{s,\beta}(N)$.

\begin{prop}
Let $G_{s,\beta}(N)$ be the graph obtained by long-range
percolation on the cycle of length $N$. Then there exists $c =
c(s,\beta) > 0$ such that
$$ \lim_{N \to \infty} \Pr \Bracket{ \tau(G_{s,\beta}(N)) \leq
\log^{c} (N) \cdot N^{s-1} } = 1 . $$ %
\end{prop}

\begin{proof}
Set $G = G_{s,\beta}(N)$.  With hindsight, set $L = \left\lceil
N^{s-1}  \xi(N)  \right\rceil$ for $\xi(N) = \alpha \log(N) / 2^s \beta$,
and $\alpha>0$ some constant to be determined below.
Let $\ell = N \pmod L$, and set $k =
\frac{N-\ell}{L}$. Divide the cycle into $k$ intervals,
$S_1,\ldots, S_k$, so that $S_1, \ldots, S_{k-1}$ are of length
$L$, and $S_k$ is of length $L + \ell \leq 2L$.
Since $s<2$, we can take $N$ large enough so that 
$2^s \beta L^2 \leq N^s$. 

Let $1 \leq i \neq j \leq k$.  
Let $\Ee(i,j)$ be the event that there exist $x \in S_i$
and $y \in S_j$ such that $(x,y) \in E(G)$.

Let $\Gamma$ be the graph obtained from $G$ by contracting 
each of the intervals $S_1,\ldots,S_k$ to a vertex.
That is $V(\Gamma) = \set{1,2,\ldots,k}$ and 
$(i,j) \in E(\Gamma)$ if $\Ee(i,j)$ occurs.

We bound from below the probability that $(i,j) \in E(\Gamma)$.
\begin{align*} 
\Pr [ (i,j) \not\in E(\Gamma) ] & = \Pr [ \textrm{not } \Ee(i,j) ] 
= \prod_{\substack{x \in S_i \\ y \in S_j } } \Pr [ (x,y) \not\in E(G) ] \\
& \leq \prod_{\substack{x \in S_i \\ y \in S_j } } \exp ( - 2^s \beta N^{-s} ) 
\leq  \exp ( - 2^s \beta L^2 N^{-s} ) .
\end{align*}
Using the inequality $1-e^{-\zeta} \geq \zeta/e$ (valid for $\zeta \in [0,1]$),
we get that 
for all $1 \leq i \neq j \leq k$,
\begin{align*} 
\Pr [ (i,j) \in E(\Gamma) ] & \geq 1 - \exp ( - 2^s \beta L^2 N^{-s} ) \\
& \geq \frac{2^s \beta}{e} \cdot \frac{L^2}{N} . 
\end{align*}
Since $k \geq \frac{N}{L}-1$, we have that for large enough $N$ (depending on $s,\beta$),
\begin{align} \label{eqn:prob. of (i,j) in Gamma}
\Pr [ (i,j) \in E(\Gamma) ] \geq \frac{\alpha}{2e} \cdot \frac{\log k}{k} , 
\end{align}
for all $1 \leq i \neq j \leq k$.

Let $\Gamma'$ be the Erdos-Renyi random graph on $k$ vertices, with edge probability 
$p = \frac{\alpha}{2e} \cdot \frac{\log k}{k}$.
That is, $V(\Gamma') = \set{1,2,\ldots,k}$ and
$(i,j) \in E(\Gamma')$ and $(j,i) \in E(\Gamma')$ with probability $p$, 
all edges $\set{i,j}$ independently.

By \eqref{eqn:prob. of (i,j) in Gamma}, we can couple $G$ and $\Gamma'$
so that $\Gamma'$ will be a subgraph of $\Gamma$.

$\deg_{\Gamma'}(j)$ has the binomiaml distribution with parameters $k,p$.
Thus, a quick calculation shows that 
there exists a constant $c_1 = c_1(\alpha) > 0$ such that 
with probability tending to $1$,
\begin{align} \label{eqn:Gamma' max degree}
\max_{1 \leq j \leq k} \deg_{\Gamma'}(j) \leq c_1 \log N . 
\end{align}
%the maximal degree in $\Gamma'$ is bounded by $c_1 \log N$.
Furthermore,
in \cite{BHKL_isoper} it is shown that for large enough
$\alpha$, there exists a constant $c_2 > 0$ such that 
with probability tending to $1$,
\begin{align} \label{eqn:mixing Gamma'}
\tau(\Gamma') \leq c_2 \log k \leq c_2 \log N .
\end{align}

We now derive an upper bound on the mixing time of $G$, 
by constructing a flow on $G$, using the optimal flow for $\Gamma'$.

Let $\pi_G,\pi_{\Gamma'}$ denote the stationary distribution of 
$G, \Gamma'$ respectively.
Let $\P(G),\P(\Gamma')$ be the set of simple paths in $G,\Gamma'$ respectively.
Let $\P(x,y;G), \P(i,j;\Gamma')$ be the set of simple paths in $G,\Gamma'$ respectively,
from $x$ to $y$, $i$ to $j$, respectively.
For a path $p$ let $p^+$ be the starting vertex of $p$, and let $p^-$ be the ending 
vertex of $p$ (specifically for edges $e = (e^+,e^-)$).

For $(i,j) \in E(\Gamma)$ let $e(i,j)$ be a specific edge
such that $e(i,j) = (x,y) \in E(G)$ and $x \in S_i$ and $y \in S_j$
(by definition there always exists at least one such edge).
For $1 \leq j \leq k$ let $G_j$ be the induced subgraph on $S_j$.
For $x,y \in S_j$ let $p(x,y)$ be a path in $G_j$ that realizes the 
distance between $x$ and $y$ in $G_j$ (a geodesic).  If $x=y$ let $p(x,x)$
be the empty path.

For $q \in \P(i,j;\Gamma')$, and $x \in V_i, y \in V_j$, 
define 
$p(q,x,y) \in \P(x,y;G)$ by interpolating $q$ using the specified edges $e(i,j)$
and geodesics $p(x,y)$;
that is if $q = e_1 e_2 \cdots e_{|q|}$,
then
$$ p(q,x,y) = p(x,e_1^+) e(e_1^+,e_1^-) p(e_1^-,e_2^+) e(e_2^+,e_2^-)
\cdots e(e_{|q|}^+,e_{|q|}^-)  p(e_{|q|}^-,y) . $$
Setting 
$\Delta = \max_j \mathrm{diam}(G_j)$
we get that $|p(q,x,y)| \leq (\Delta+1) |q|$.

Let $f^*$ be the optimal flow on $\Gamma'$.  As mentioned above,
by Theorem 8 of \cite{Sinclair}, using also \eqref{eqn:mixing Gamma'},
there exists a constant $c_3>0$ such that
with probability tending to $1$,
\begin{align} \label{eqn:flow Gamma'}
\forall \ (i,j) \in E(\Gamma') \qquad
|E(\Gamma')| \sum_{\substack{ q \in \P(\Gamma') \\ q \ni (i,j) } }
f^*(q) |q| & \leq 16 \tau(\Gamma') 
\leq c_3 \log N . 
\end{align}

We now define a flow $f$ on $G$ using $f^*$.
Let $x,y \in V(G)$, and let $i,j$ be such that $x \in S_i$ and $y \in S_j$.

If $i=j$ set $f(p) = \pi_G(x) \pi_G(y)$ if $p = p(x,y)$ and $0$ otherwise.

If $i \neq j$, then for any $q \in \P(i,j;\Gamma')$ set 
$$ f(p(q,x,y)) = \frac{f^*(q)}{\pi_{\Gamma'}(i) \pi_{\Gamma'}(j)} 
\cdot \pi_G(x) \pi_G(y) , $$
and $0$ otherwise.

We calculate the conestion of the flow $f$.
Let $(x,y) \in E(G)$, and let $i,j$ be such that $x \in S_i$ and $y \in S_j$.

{\bf Case 1:} $i \neq j$.
In this case, any path $p$ that contains the edge $(x,y)$,
such that $f(p) > 0$, 
must be of the form $p = p(q,z,w)$ for some $q \in \P(\Gamma')$ that contains $(i,j)$.
Thus, using \eqref{eqn:Gamma' max degree} and \eqref{eqn:flow Gamma'},
\begin{align} \label{eqn:long edge}
\sum_{\substack{p \in \P(G) \\ p \ni (x,y)} } f(p) |p| & \leq 
\sum_{\substack{q \in \P(\Gamma') \\ q \ni (i,j) } }
\sum_{\substack{ z \in S_{q^+} \\ w \in S_{q^-} } }
\frac{\pi_G(z) \pi_G(w)}{\pi_{\Gamma'}(q^+) \pi_{\Gamma'}(q^-) } \cdot 
f^*(q) |q| (\Delta+1) \nonumber \\
& \leq (\Delta+1) \cdot (\max_j \pi_G(S_j))^2 \cdot |E(\Gamma')|^2 \cdot 
\sum_{\substack{q \in \P(\Gamma') \\ q \ni (i,j) } } f^*(q) |q| \nonumber \\
& \leq \frac{1}{|E(G)|} (\Delta+1) \cdot (\max_x \deg_G(x))^2 \cdot \frac{L^2 k}{N}
 \cdot c_4 \log^{c_5} N \nonumber \\
& \leq \frac{1}{|E(G)|} \cdot c_4 \log^{c_5} N \cdot 
(\Delta+1) \cdot (\max_x \deg_G(x))^2 \cdot L ,
\end{align}
where $c_4,c_5>0$ are constants (independent of $N$).

{\bf Case 2:} $i=j$.
I this case, any path $p$ that contains the edge $(x,y)$,
such that $f(p) > 0$, 
is one of the follwing:
Either it is of the form $p = p(q,z,w)$ for some $q \in \P(\Gamma')$ that
contains the vertex $i$, or it is of the form $p = p(z,w)$ for some $z,w \in S_i$.
Any path $q \in \P(\Gamma')$ that contains the vertex $i$ must contain some edge
$(i,j) \in E(\Gamma')$.  Thus, using \eqref{eqn:Gamma' max degree}
and \eqref{eqn:long edge},
\begin{align} \label{eqn:short edge}
\sum_{\substack{p \in \P(G) \\ p \ni (x,y)} } f(p) |p| 
& \leq \sum_{z,w \in S_i} f(p(z,w)) |p(z,w)| + \sum_{j:(i,j) \in E(\Gamma')}
\sum_{\substack{q \in \P(\Gamma') \\ q \ni (i,j) } } 
\sum_{\substack{ z \in S_{q^+} \\ w \in S_{q^-} } }
f(p(q,z,w)) |p(q,z,w)|
\nonumber \\
& \leq 
\sum_{z,w \in S_i} \pi_G(z) \pi_G(w) \Delta + \deg_{\Gamma'}(i) 
\cdot \frac{1}{|E(G)|} \cdot c_4 \log^{c_5} N \cdot 
(\Delta+1) \cdot (\max_x \deg_G(x))^2 \cdot L \nonumber \\
& \leq \frac{1}{|E(G)|} \cdot c_6 \log^{c_7} N \cdot 
(\Delta+1) \cdot (\max_x \deg_G(x))^2 \cdot L , 
\end{align}
where $c_6,c_7>0$ are constants (independent of $N$).

By our choice of $L$, and by \eqref{eqn:mixing}, it suffices to show that
there exists constants $c_8,c_9,c_{10},c_{11} >0$ (that may depend on $s,\beta$) such that
$\Delta \leq c_8 \log^{c_9} N$ and $\max_x \deg_G(x) \leq c_{10} \log^{c_{11}} N$,
with probability tending to $1$.

In \cite{LRP}, the following is shown:
There exist $\delta = \delta(s,\beta) > 0$
and $n_0 = n_0(s,\beta) > 0$ such that 
for all $n > n_0$,
$$ \Pr [ \mathrm{diam} (\textrm{long-range
percolation on the interval of length } n) > \log^\delta(n) ]
< \frac{1}{n^2} . $$
% where $\delta =
% \frac{\log 2}{\log(2/s)} + o(1)$.
Thus, since there are $k \leq N/L$ such intervals, each of length at least $L$,
a union bound gives that for large enough $N$,
\begin{align} \label{eqn:diameter}
\Delta = \max_j \mathrm{diam}(G_j) \leq \log^\delta(N) ,
\end{align}
with probability at least $1-O(N^{4-3s})$, which tends to $1$.

Next, we show that with probability tending to $1$, the maximal
degree in $G$ is bounded by $2 \log(N)$.  This follows from the
following considerations:

Fix a vertex $x$ in $G$.  We can
write $\deg_G(x) = 2 + \sum_{y \neq x} Z_{xy}$, where $Z_{xy}$ is the
indicator function of the event that $x$ and $y$ are connected by
an edge not in the cycle.  The random variables $Z_{xy}$ are all
independent.  Note that $\mu_{xy} = \E \br{ Z_{xy} } = 1- e^{-
\beta \norm{x-y}^{-s} }$, and for any $t>0$, $\E \br{ \exp \sr{ t
Z_{xy} } } = 1 + \mu_{xy} \sr{ e^t -1} \leq \exp \sr{ (e^t-1)
\mu_{xy} }$. So we can calculate using Markov's inequality, for
all $\lambda > 0$,
\begin{eqnarray*}
    \Pr \br{ \deg_G(x) > 2 + \lambda } & \leq & 
    \exp \sr{ -\lambda } \prod_{y \neq x} \exp \sr{ (e-1)
    \mu_{xy} } \\
%     & \leq & \exp \sr{ - \lambda + (e-1) \beta \sum_{y \neq x}
%     \norm{ x-y}^{-s} } \\
%     & \leq & \exp \sr{ - \lambda + (e-1) 2 \beta
%     \sum_{k=1}^{N/2} k^{-s} } \\
    & \leq & \exp \sr{ -\lambda + c(s,\beta) N^{1-s} } .
\end{eqnarray*}
Thus, taking $\lambda = 2 \log(N) -2$ and using a union bound, the
probability that there exists a vertex $x$ with $\deg_G(x) > 2 \log(N)$
is bounded by $N \frac{1+o(1)}{N^2} \leq \frac{1+o(1)}{N}$.  Thus,
with probability tending to $1$ we have that
\begin{align} \label{eqn:degree G}
\max_{x \in V(G)} \deg_G(x) \leq 2 \log (N) . 
\end{align}

Finally, combining 
Theorem 5' of \cite{Sinclair} with the flow $f$,
\eqref{eqn:mixing}, \eqref{eqn:long edge}, \eqref{eqn:short edge},
\eqref{eqn:diameter} and \eqref{eqn:degree G}, and by our choice of $L$, 
we conclude that 
there exist constants $c,c'>0$, (independent of $N$, but perhaps depending on $s,\beta$),
such that with probability tending to $1$,
$$ \tau(G) \leq c \log{c'}(N) \cdot N^{s-1} . $$ 
\end{proof}

\section{Lower Bound}

In Corollary 5.2 of \cite{LRP}, it is shown that $\tau(G) \geq
\Omega \sr{ N^{s-1} }$ with probability tending to $1$. This
miraculously matches our upper bound up to poly-logarithmic
factors. In \cite{LRP} there is a typo in the parameters, so for
completeness we provide the proof here.

\begin{prop}
Let $G_{s,\beta}(N)$ be the graph obtained by long-range
percolation on the cycle of length $N$. Then
$$ \lim_{N \to \infty} \Pr \Bracket{ \tau(G_{s,\beta}(N)) \leq c N^{s-1}  } = 0 , $$ %
where $c = c(s,\beta)$ is a constant independent of $N$.
\end{prop}

\begin{proof}
Let $G = G_{s,\beta}(N)$.  It is well known (see e.g.
\cite{Aldous, DS, Sinclair}) that it is enough to bound from above
the Cheeger constant of the graph $G$, which is defined as
$$ {\cal C}(G) = \min_{\substack{ \emptyset \neq A \subset V(G) \\
\abs{A} \leq N/2}} \frac{\abs{\partial A}}{\abs{A}} \qquad
\partial A = \setb{ \set{x,y} \in E(G) }{ x \in A, y \not\in A } $$

The natural subset to choose is $A = \set{1,2,\ldots, N/2}$ (any
arc of length $N/2$ will suffice).  For $x \in A$ and $y \not\in
A$ let $Z_{xy}$ be the indicator function of the event that $x$
and $y$ are connected by an edge not in the cycle.  Then, $Z_{xy}$
are all independent.  Set $\mu_{xy} = \E [ Z_{xy} ] = 1- \exp
\Soger{ - \beta \norm{x-y}^{-s} }$. For any $t > 0$,
$$ \E \Bracket{ \exp (t Z_{xy}) } = 1 + \mu_{xy} (e^t-1) \leq \exp
(\mu_{xy} (e^t -1) ) . $$ %
We have that $\abs{\partial A} = 2 + \sum_{x \in A} \sum_{y
\not\in A} Z_{xy}$, and that
\begin{eqnarray*}
    \sum_{x \in A} \sum_{y \not\in A} \mu_{xy} & \leq &
    \sum_{x \in A} \sum_{y \neq x} \mu_{xy} \\
    & \leq & \sum_{x \in A} c(s,\beta) N^{1-s} \leq c(s,\beta) N^{2-s} .
\end{eqnarray*}

Thus, there exists $c = c(s,\beta)$ such that for any $t>0$ and
any $\lambda > 0$,
\begin{eqnarray*}
    \Pr \br{ \abs{\p A} > 2 + \lambda N^{2-s} } & \leq &
    \frac{ \E \br{ \exp \sr{ t (\abs{\p A} - 2) } } }{ \exp \sr{
    t \lambda N^{2-s} } } \\
    & = & \exp \sr{ - t \lambda N^{2-s} }
    \prod_{x \in A} \prod_{y \not\in A } \E \Bracket{ \exp (t Z_{xy})
    }  \\
    & = & \exp \sr{ (e^t-1) \sum_{x \in A} \sum_{y \not\in A}
    \mu_{xy} } \exp \sr{ - t \lambda N^{2-s} } \\
    & \leq & \exp \sr{ N^{2-s} \cdot \sr{ c(s,\beta) (e^t-1) - t
    \lambda } } .
\end{eqnarray*}
Choosing $t$ small enough, we get that for some fixed $c =
c(s,\beta)$ we have that
$$ \lim_{N \to \infty} \Pr \br{ \abs{\p A} > c N^{2-s} } = 0 , $$
which implies that with probability tending to $1$ as $N$ tends to
infinity, ${\cal C}(G) \leq c N^{1-s}$ for some (possibly
different) $c = c(s,\beta)$ independent of $N$.

This gives a bound on the mixing time (see \cite{Sinclair}):
$$ \tau(G) \geq \frac{1-\log 2}{2} \Soger{ \frac{1}{2 {\cal C}(G)}
-1} \geq c N^{s-1} , $$ %
with probability tending to $1$, and $c = c(s,\beta)$ independent
of $N$.
\end{proof}

\section{A Phase Transition}  \label{phase}

In the previous sections we have shown that the mixing time of
$G_{s,\beta}$ is $N^{s-1}$ (disregarding poly-logarithmic factors)
for $1<s<2$.   When $s$ tends to $2$, this quantity tends to $N$.
We will show that a phase transition occurs at $s=2$, meaning that
for $s>2$ the mixing time will ``jump'' to $N^2$.

\begin{prop}
Let $G = G_{s,\beta}(N)$ for $s>2$.  Then the mixing time of $G$
satisfies
$$ \lim_{N \to \infty} \Pr \Bracket{ \tau(G) \geq c N^2 } = 1 $$ %
for some constant $c = c(s,\beta)$, independent of $N$.
\end{prop}

\begin{proof}
For simplicity we assume that $N$ is divisible by $8$. For other
$N$ the proof is similar.

Set
$$ A = \set{1, 2, \ldots, \frac{N}{2} } \quad B = \set{ \frac{N}{2} + 1 , \ldots, \frac{3 N}{4} }
\quad C = \set{ \frac{3 N}{4} + 1, \ldots, N} . $$ %
Also, for $i=1,2,\ldots,8$ set $K_i = \set{ (i-1) \frac{N}{8} +1,
\ldots, i \frac{N}{8} }$.

By the proof of Theorem 3.1 (A) in \cite{LRP}, there exists $c_1 =
c_1(s,\beta)$ such that with probability $1-o(1)$, we have that
all sets $K_1,\ldots,K_8$ each contain at least $c_1 N$ vertices
of degree $2$.

Further, by rotating the cycle, without loss of generality we can
assume that $\pi(A) \geq \pi(B \cup C)$, and that $\pi(B) \geq
\pi(C)$.  This implies that $\pi(A \cup B) \geq \frac{3}{4}$.

Fix a vertex $x \not\in A \cup B$, and let $\sr{ S_t \ ; \ t \geq
0 }$ be a simple random walk on $G$ starting at $S_0 = x$. Let $T$
be the hitting time of the set $A \cup B$. Note that at any time
$t \geq \tau(G)$, we have that
$$ \frac{3}{4} - \Pr \br{ S_t \in A \cup B } \leq \sum_{y \in
A \cup B} \abs{ \pi(y) - P^t (x,y) } \leq 2 \Delta_x (t) \leq
\frac{2}{e} . $$ %
So we conclude that for any $x \not\in A \cup B$ and $t \geq
\tau(G)$ we have that
$$ \Pr_x \br{ T \leq t } \geq \Pr_x \br{ S_t \in A \cup B } > 0.01
. $$ %
This implies that for any $x \not\in A \cup B$ and for any real
$s$, $\Pr_x \br{ T > s } \leq 0.99^{(s/\tau -1)}$. Thus, there
exists $c_2 > 0$ independent of $N$, such that
$$ \E_x \br{ T } = \sum_{t=0}^{\infty} \Pr_x \br{ T > t } \leq c_2 \tau . $$

Set $u = \frac{3}{8} N$.  Recall that there are at least $c_1 N$
vertices of degree $2$ separating $u$ from $A \cup B$ (on each
side of $u$).  We will show that this implies that $\E_u \br{T}
\geq c_3 N^2$ for some $c_3 = c_3(s,\beta)$ independent of $N$.

We use the language of electrical networks, see \cite{DoyleSnell,
LyonsPeres} for background.  We remark that for the reader not
familiar with these notions, one can use the Varopoulos-Carne
bounds (see e.g. \cite{LyonsPeres}) to show that a linear diameter
implies that the mixing time is at least $\frac{c N^2}{\log N}$.

We can write $T = \sum_{x \not\in A \cup B} V_x$ where $V_x$ is
the number of visits to the vertex $x$, up to time $T$.  Ground
the set $A \cup B$ (so that its voltage is $0$), and set a
potential to $u$ so that there is a unit current flowing into $u$.
We get that for any $x$, we have the identity $\E_u \br{ V_x} =
v(x) d(x)$, where $v(x)$ is the voltage at $x$, and $d(x)$ is the
degree of $x$ (this follows from noting that $\E_u \br{V_x} /
d(x)$ is harmonic).

Let $x_1,x_2, \ldots, x_{c_1 N}$ be the vertices of degree $2$ on
the side of $u$ with at least $1/2$ the current (without loss of
generality say in the interval $\set{ \frac{3}{8}N + 1, \ldots,
N}$). Fix $1 \leq i \leq c_1 N$. Since $x_i$ is a cut point, the
current flowing into and out of $x_i$ is at least $1/2$. Thus, the
voltage at $x_i$, $v(x_i)$, is at least $1/2$ the resistance
between $x_i$ and $A \cup B$. This resistance is bounded from
below by the number of cut \emph{edges} (which are resistors of
resistance $1$ connected serially), which in turn is bounded by
the number of vertices of degree $2$ between $x_i$ and the set $A
\cup B$. That is, $v(x_i) \geq i/2$. Thus,
$$ \E_u \br{T} \geq \sum_{\substack{x \not\in A \cup B \\ d(x) = 2
} } \E_u \br{ V_x} = \sum_{\substack{x \not\in A \cup B \\ d(x)
= 2 } } d(x) v(x) \geq \sum_{i=1}^{c_1 N} i  = c_3 N^2 , $$ %
for some $c_3 = c_3(s,\beta)$ independent of $N$.

We have shown that with probability $1-o(1)$ there are a linear
number of vertices of degree $2$ separating a vertex $u$ from a
set of high weight under the stationary distribution.  Thus,
$\tau(G) \geq c_4 \E_u \br{T} \geq c_5 N^2$, where $c_4,c_5$
depend only on $s$ and $\beta$.
\end{proof}

{\bf Acknowledgement.} We wish to thank Ori Gurel-Gurevich for
useful conversations.

% -----------------------------------------------------------------------------------
% ---------------------------- BIBLIOGRAPHY -----------------------------------------
% -----------------------------------------------------------------------------------

\end{document}